\documentclass[12pt]{article}
\usepackage{amssymb,amsmath,amsthm}
\begin{document}

\begin{center}
\textbf{ON GENERALIZATION OF THE FREUDENTHAL'S THEOREM FOR COMPACT
IRREDUCIBLE STANDARD POLYHEDRIC REPRESENTATION FOR
SUPERPARACOMPACT COMPLETE METRIZABLE SPACES}\\

[D.K.Musaev], D.I.Jumaev
\end{center}
\begin{abstract}
In this paper for superparacompact complete metrizable spaces the
Freudenthal's theorem for compact irreducible standard polyhedric
representation is generalized. Furthermore, for superparacompact
metric spaces are reinforced: 1) the Morita's theorem about
universality of the product $Q^\infty\times B(\tau)$ of Hilbert
cube $Q^\infty$ to generalized Baire space $B(\tau)$ of the weight
$\tau$ in the space of all strongly metrizable spaces of weight
$\le \tau$; 2) the Nagata's theorem about universality of the
product $\Phi^n\times B(\tau)$ of universal $n$- dimensional
compact $\Phi^n$ to $B(\tau)$  in the space of all strongly
metrizable spaces $\le\tau$ and dimension $dimX\le n.$
\end{abstract}

In what follows, by space we mean topological spaces, by compacts --
metrizable bicompacts, mapping is used to mean - continuous mapping
of spaces. Furthermore, in this paper as polyhedron we mean spaces
(generally speaking, infinity) simplicial complex (see. [1], Chapter
3 \S2) in metrizable topology.

We give main definitions and some necessary concepts for this
paper.

\textbf{Definition 1.} [1] a) A system $\omega$ of subsets of the
set $X$ is called star countable (finite), if every element of the
system $\omega$ is intersected at most countable (finite) number
of elements  of this system; b) finite sequence of subsets
$M_0,...,M_s$, of the set $X$ is called chain, connecting sets
$M_0$ and $M_s$, if $M_{i-1}\cap M_i\neq \emptyset$ for all values
$i=1,...,s$; c) a system $\omega$ of subsets of the set $X$ is
called enchained set, if for all sets $M$ and $M'$ of this system
there exists such chain of elements of the system $\omega$ that
the first element of the chain is the set $M$, and the last is the
set $M'$; maximal enchained subsystems of the system $\omega$ is
called components of overlapping (or components) of the system
$\omega$.

Known that [see [1], Chapter 1, \S 6] components of the star
countable system $\omega$ are countable.

{\bf Definition 2} [2] (a) A star-finite open covering of a space
is called finite-component if all components of overlapping are
finite;

(b) A space is called superparacompact if, given an open covering of
the space, a finite-component covering can be inscribed in this open
covering;

(c) Hausdorff superparacompact spaces are called superparacompacta.

{\bf Definition 3.} [1]. (a) A finite covering $\omega=\{ O_1,
\dots, O_s\}$ of a space $X$ is said to be \emph{irreducible}, if
no proper subcomplex $N'$ of the nerve (see [1, Chapter 3, \S 2])
$N_{\omega}$ of the covering $\omega$ is not nerve of covering
more small than $\omega$ (i.e. of covering, inscribed into
covering $\omega$);

(b) complete finite complex $K$ (see [1, Chapter 3, \S2]),
elements of which are disjoint open simplexes of given space $R^n$
is called \emph{triangulation}, lying in $R^n$.

(c)) the mapping  $f$ of the space $X$ in body of triangulation {\it
irreducible} with respect to this triangulation, if it essentially
(see [1, Chapter. 3, \S 5]) on preimage of every closed simplex of
this triangulation;

(d) Let the triangulation $N=N_{\omega}$ be geometrical
realization into $R^m$ of nerve of the covering
$\omega=\{O_1,...,O_s\}$ of the space $X$ and $e_i$ be vertex of
the nerve $N$ correspondint to a element $O_i$ of the covering
$\omega$.

The mapping $f$ of the space $X$ to polyhedron $\widetilde N$ (see
[1, Chapter 4, \S 1]) is called {\it finite} with respect to
covering $\omega,$ if preimage $f^{-1} O{e_i}$ of every star
$Oe_i$ contains in $O_i$;

(e) the space with $\sigma$-star finite base is called (see [1,
Chapter 6, \S 3]) {\it strongly metrizable}.

{\bf Definition 4.} Finite component covering $\omega$ of the space
$X$ is called {\it irreducible}, if all its components coupling are
irreducible (i.e. components $\omega_\lambda$ of covering  $\omega$
are irreducible covering of own bodies $\widetilde \omega_\lambda$).

{\bf Definition 5} [1]. {\it By product} of two systems of sets
  $\alpha=\{ A\}$ and $\beta=\{ B\}$ is called the system of sets $\gamma=\alpha\wedge\beta$, by elements of
which are all (denoted) of sets of the form $A\cap B,$ where $A\in
\alpha,$ $B\in \beta.$

\textbf{Preposition 1.} \emph{Into any open covering of
superparacompact Hausdorff space $X$ can be inscribed irreducible
finite-component open covering.}

\textbf{Proof.} Let $\omega$  be any open covering of a
superparacompact Hausdorff space $X.$ Since the space $X$
superparacompact, then, without limiting of generality of argument,
the covering $\omega$ can be counted  finite component. Since every
component $\omega_\lambda,$ $\lambda\in L$ of covering $\omega$ is
finite, and their body $\widetilde \omega_\lambda,$ $\lambda\in L$
is open-closed in the $X,$ then in the covering $\omega_\lambda$ of
the set   $\widetilde \omega_\lambda$ can be inscribed (see [1,
preposition 2, Chapter 4, \S2]) irreducible open covering
$\omega_\lambda^*.$ Then the system $\omega^*=\cup\{
\omega_\lambda^*: \lambda\in L\}$ is irreducible finite component
open covering of the space $X,$ inscribed into $\omega.$ \hfill
$\Box$

{\bf Remark 1.} a) If $\omega=\{ O_\alpha, \alpha\in A,
|A|=\tau\}$ is finite component open covering of the space
  $X,$ then the body $\widetilde{gN}_\omega$ of standard geometric realization
$gN_\omega$ in Hilbert space  $R^\tau$ of nerve $N_\omega$ of
covering  $\omega$ (standardization of realization means, that
vertices of the triangulation $gN_\omega$ are located in unit points
of the space $R^\tau$) is discrete subcompact polyhedrons
$\widetilde{gN}_{\omega_\lambda}$, being body of realization of
nerves  $N_{\omega_\lambda}$ of component $\omega_\lambda$ of
covering $\omega.$

We note, that polyhedron  $\widetilde{gN}_\omega$ is
superparacompact  (see. [3, proposition 2]). If, furthermore in
addition the covering $\omega$ has multiplicity  $\le n+1$, then
polyhedrons $\widetilde{gN}_{\omega_\lambda}$ is not more than
$n$-dimensional and thus  $dim\widetilde{gN}_\omega\le n$.

b) From the theorem about canonical mappings   (see [1, Chapter 4,
\S 4, theorem 1]) when transfer to bodies of components of finite
component covering it is easy to get following

\textbf{Preposition 2.} Let $\omega=\{ O_\alpha,\alpha\in A\}$ be
any finite component open covering  of the normal space
  $X$ with nerve  $N_\omega,$ realized in the triangulation form; it is possible to find such subcomplex   $N_\omega'$
  of the nerve $N_\omega$ and the mapping $f:X\to \widetilde N_\omega,$ which canonical with respect to  $\omega$
 in which the image $fX$ is polyhedron
$\widetilde N_\omega'\subseteq \widetilde N_\omega$ and every
principal simplex of the complex  $N_\omega'$ is covered
essentially.

\textbf{Preposition 3.} For any finite component irreducible
covering $\omega=\{ O_\alpha, \alpha\in A, |A|=\tau\}$ of the normal
space $X$ arbitrary canonical mapping of the space $X$ into body
$\widetilde{gN}_\omega$ of standard geometrical realization
$gN_\omega$ in Hilbert space  $R^\tau$ of the nerve
 $N_\omega$ of covering  $\omega$ is irreducible  with respect to
 triangulation  $gN_\omega$.

\textbf{Proof.} Let $f$ be any canonical mapping  (with respect to
covering $\omega$) of the space $X$ into body
$\widetilde{gN}_\omega$ of standard geometrical realization
$gN_\omega$ in Hilbert space  $R^\tau$ of the nerve $N_\omega$ of
covering$\omega=\{ O_\alpha, \alpha\in A, |A|=\tau\}$. Then,
according to remark 1, the body $\widetilde{gN}_\omega$ of
standard geometrical realization $gN_\omega$ into  $R^\tau$ of the
nerve  $N_\omega$ of the covering $\omega$ is discrete sum compact
polyhedrons $\widetilde{gN}_{\omega_\lambda}$, being bodies of
realization of nerves $N_{\omega_\lambda}$ of components
${\omega_\lambda}$ of covering $\omega.$ Suppose $f_\lambda=f:
\widetilde \omega_\lambda\to \widetilde{gN}_{\omega_\lambda}$ for
any $\lambda\in L.$ Clearly, that every mapping $f_\lambda,$
$\lambda\in L,$ is canonical  (with respect to the covering
$\omega_\lambda$).

Since the covering $\omega$ is nonreducible, then all its components
$\omega_\lambda,$ $\lambda\in L,$ are irreducible, according to
definition 3, and therefore each canonical mapping $f_\lambda:
\widetilde \omega_\lambda\to \widetilde{gN}_{\omega_\lambda}$ is
(see [1, Chapter 4, \S3]) irreducible. Then canonical mapping
$f:X\to \widetilde{gN}_\omega$, as combination [4] of irreducible
mappings $\{ f_\lambda, \lambda\in L\}$, is irreducible with respect
to triangulation $gN_\omega.$ \hfill $\Box$

Later in this work as triangulation in the Hilbert space $R^\tau$
is meant either standard geometrical realization  $gN_\omega$ of
nerve  $N_\omega$ of finite component covering
 $\omega$ of normal space, or such its subdivision
 (see [1, Chapter 3, \S 2, section 5]) $(gN_\omega)^*$, which for each component  $\omega_\lambda$
of the covering $\omega$ coincide with some (multiple, and also
multiplicity  depends on component $\omega_\lambda$) barycentric
subdivision (see [1, Chapter 3, \S 2, section 6])of the
triangulation $gN_{\omega_\lambda}.$

{\bf Remark 2.} a) Let $\omega=\{ O_\alpha,\alpha\in A, |A|\le
\tau\}$ be finite component  $(n+1)$-multiplicity covering of
normal space $X,$ $\widetilde{gN}_{\omega_\lambda}$ is body of
standard geometrical realization $gN_\omega$ in the Hilbert space
$R^\tau$ of the nerve  $N_\omega$ of the covering $\omega$ and
$\varepsilon>0.$ We take such natural $s,$ that $ \left(
\frac{n}{n+1}\right)^s\sqrt 2 < \varepsilon.$ Then in virtue of
isometrically of all $k$-dimensional simplexes of the
triangulation $gN_\omega$ and relation  $ \left(
\frac{k}{k+1}\right)^s \sqrt 2\le  \left( \frac{n}{n+1}\right)^s
\sqrt 2,$ $k=1,2,\dots,n,$ follows, that all simplexes of
subdivision  $\left( gN_\omega\right)^*$, being $s$-multiplicity
barycentric subdivision of the triangulation $gN_{\omega}$, have
diameter $<\varepsilon.$

b) Let $\omega_1$  be finite component open covering of the normal
space   $X,$ $gN_{\omega_1}$ is standard geometrical realization
of the nerve $N_{\omega_1}$ of the covering   $\omega_1$ in the
Hilbert space  $R^\tau$ and $f_1$ is canonical with respect to the
covering $\omega_1$ mapping of the space $X$ into
$\widetilde{gN}_{\omega_1}$. Let $\left( gN_{\omega_1}\right)^*$
be triangulation of the polyhedron $\widetilde{gN}_{\omega_1}$,
being subdivision of the triangulation $gN_{\omega_1},$ and
covering   $\omega_2'$ consists on preimages in the mapping
  $f_1$ of main stars  (see [1, Chapter 3, \S 2, section 3])
of the triangulation   $\left( gN_{\omega_1}\right)^*$. Suppose
also that finite component covering $\omega_2$ of  the space
  $X$ inscribed into covering $\omega'_2,$ $gN_{\omega_2}$
is standard geometrical realization of  the nerve $N_{\omega_2}$
and $f_2$ is canonical with respect to   $\omega_2$ mapping of the
space $X$ into $\widetilde{gN}_{\omega_2}$.

Then any generated refinement $\omega_2$ into $\omega_2'$ (see [1,
Chapter 3, \S 2, section 4]) simplicial with respect to
triangulation   $gN_{\omega_2}$ and
$\left({gN}_{\omega_1}\right)^*$ the mapping $\pi:
\widetilde{gN}_{\omega_2} \to \widetilde{gN}_{\omega_1}$ is
obtained, according to the lemma (see [1, Chapter 3, \S1]) about
descent with respect to triangulation
$\left({gN}_{\omega_1}\right)^*$ from mapping
  $f_1$ (i.e. support arbitrary point   $\pi f_2(x)$ is face of
  support of the point   $f_1(x)$ in the triangulation
$\left({gN}_{\omega_1}\right)^*$).

The proof implies  from the case of compact polyhedrons  (see [1,
Chapter 3, \S 1]) we turn on to bodies of component of the
covering  $\omega_2.$

\textbf{Theorem 1.} \emph{Any $n$-dimensional complete metric
superparacompact space $X$ is  limit of inverse sequence $S=\left\{
\widetilde K_i, \pi_i^{i+1}\right\},$ $i=1,2,\dots,$ from
$n$-dimensional polyhedrons $\widetilde K_i,$ being bodies of
standard triangulation $K_i$  decomposing to discrete sum of compact
polyhedrons; in addition projections $\pi_i^{i+1}$ are simplicial
with respect to $K_{i+1}$ and some triangulation $K_i^*$ of the
polyhedron $\widetilde K_i$, being subdivision of the triangulation
$K_i.$ Every projection $\pi_i: X\to \widetilde K_i$ is irreducible
with respect to triangulation $K_i$, $i=1, 2, \dots$.}

\textbf{Proof.} We construct searching inverse sequence by
induction. Let $\gamma_i,$ $i=1, 2, \dots,$ be $1/{2^i}-$ open
covering of the space  $X.$ Since $dimX=n,$ then there exists such
open covering  $\eta$ of the space $X,$ that any inscribed covering
into it has multiplicity  $\ge n+1.$ By virtue of preposition 1, to
the covering $\{\gamma_1\wedge \eta\}$ we inscribe irreducible
finite component open covering  $\omega_1$ of the space $X.$ Nerve
of the covering $\omega_1$ we denote by $N_1,$ and as $K_1$ we
denote standard geometrical realization $N_1$ in Hilbert space
$R^\tau$. According to preposition 2, there exists finite with
respect to $\omega_1$ mapping $f_1$ of the space $X$ into polyhedron
$\widetilde K_1.$ Because the covering
 $\omega_1$ is irreducible, then, according to preposition3, the mapping $f_1$ is
 irreducible mapping
with respect to triangulation $K_1$ and, so, will be mapping onto
 $\widetilde K_1.$ The covering $\omega_1$ inscribed into covering $\{ \gamma_1\wedge\eta\}$ of the space $X,$
thus the covering  $\omega_1$ has multiplicity $n+1,$ and
$dim\widetilde K_1 =n.$ As the covering $\omega_1$ is finite
component polyhedron  $\widetilde K_1$ is discrete sum of compact
polyhedrons. We consider covering $\varphi_1,$ consisting of
preimages main stars of the triangulation $K_1^*$ in the mapping
 $f_1,$ where $K_1^*$ is such subdivision of the triangulation $K_1,$ that its mesh $< 1/2^2$ (see section a)
of remark 2).

Into covering $\{\varphi_1\wedge \eta \wedge \gamma_2\}$ we inscribe
irreducible finite component open covering
  $\omega_2$ of the space $X.$ According to preposition 2, there is canonical with respect to
$\omega_2$ mapping  $f_2$  of the space $X$ into polyhedron
$\widetilde K_2,$ where $K_2$ is standard geometrical realization of
the nerve $N_{\omega_2}$ of the covering $\omega_2$ into $R^\tau.$
By that reason, that given above, canonical with respect to
$\omega_2$ mapping  $f_2$ of the space $X$ into polyhedron
$\widetilde K_2$ is irreducible with respect to triangulation
 $K_2;$ the covering
$\omega_2$ has multiplicity  $n+1;$ the polyhedron $\widetilde K_2
$ is discrete sum of compact polyhedrons and $dim\widetilde
K_2=n$. We take some generated with inscribed  $\omega_2$ in $\{
\varphi_1\wedge\eta\wedge\gamma_2\}$ simplicial with respect to
the triangulation $K_2$ and $K_1^*$ mapping $\pi_1^2: \widetilde
K_2\to \widetilde K_1$. Then, according to section b) of the
remark  2, the mapping $\pi_1^2f_2$ is descent of the mapping
$f_1$ with respect to triangulation $K_1^*.$ Therefore $d\left(
f_1, \pi_1^2 f_2\right) < 1/2^2.$

Suppose, that for all  $i<m$ we constructed: a) $n$-dimensional
polyhedrons $\widetilde K_i,$ being bodies of standard geometrical
realizations in $R^\tau$ of nerves  $N_{\omega_i}$ of irreducible
finite component coverings $\omega_i$ of the space $X,$ inscribed
into coverings  $\{\eta\wedge \gamma_i\},$ $i=1, 2, \dots$; b)
canonical with respect to  $\omega_i$ mappings $f_i:X\to\widetilde
K_i,$ being irreducible mappings with respect to triangulations
$K_i;$ c) mappings $\pi_{i-1}^i:\widetilde K_i\to \widetilde
K_{i-1},$ $2<i<m,$ which simplicial with respect to triangulation
 $K_i$ and some triangulation  $K^*_{i-1}$ of polyhedron  $\widetilde K_{i-1}$,
 being subdivision of triangulation $K_{i-1}$; in this connection the mapping
$\pi^i_{i-1}f_i$ is obtained from $f_{i-1}$ by descent with
respect to  $K^*_{i-1}$; d) mappings
$\pi_j^i=\pi_j^{j+1}\ldots\pi_{i-1}^i,$ $\pi_i^i,$ $j<i,$ satisfy
inequalities $d\left( \pi_j^{i-1} f_{i-1}, \pi_j^i f_i\right) <
1/2^i.$

Assume now $i=m.$ According to the remark  1, the polyhedron
$\widetilde K_{m-1}$ is discrete sum of compact polyhedrons
$\widetilde K_{m-1}^\beta$, $\beta\in L,$ being bodies of standard
realizations  $K^\beta_{m-1}$ into $R^\tau$ of nerves of components
of the covering $\omega_{m-1}.$ In triangulation $K^\beta_{m-1},$
$\beta\in L,$ there exists such barycentric subdivision $\left(
K^\beta_{m-1}\right)^{s(\beta)}$, that all simlexes of the
triangulation $\left( K^\beta_{m-1}\right)^{s(\beta)}$ and their
images into polyhedrons $\widetilde K_j$ in the mapping
$\pi_j^{m-1},$ $j\le i\le m-2,$, have diameters  $<1/2^m.$ Suppose
$\left( K_{m-1}\right)^*$ coinciding with  $\left( K^\beta_{m-1}
\right)^{s(\beta)}$ on $\left(\widetilde {K^\beta_{m-1}}\right)$.
Clearly, that all simplexes of the triangulation $\left(
K_{m-1}\right)^*$ and their images into polyhedrons $\widetilde K_j$
in mappings $\pi_j^{m-1},$ $j\le i\le m-2,$, have diameters
$<1/2^m.$ Into the covering $\{ \varphi_{m-1}\wedge\eta\wedge
\gamma_m\},$ where $\varphi_{m-1}$ consists on preimages of main
stars of the triangulation $K^*_{m-1}$, in the mapping $f_{m-1}$,
according to preposition 1, we inscribe irreducible finite component
open covering  $\omega_m$ of the space  $X.$ There exists canonical
with respect to $\omega_m$ mapping of the space $X$ into polyhedron
$\widetilde K_m,$ where $K_m$ is standard geometrical realization of
the nerve $N_{\omega_m}$ of the covering $\omega_m$ into $R^\tau.$
As before, canonical with respect to  $\omega_m$ mapping $f_m$ of
the space  $X$ into polyhedron $\widetilde K_m$ is irreducible with
respect to triangulation   $K_m$ (and, so, will be mapping onto
$\widetilde K_m$); the covering $\omega_m$ have multiplicity $n+1;$
polyhedron  $\widetilde K_m$ is discrete sum of compact polyhedrons
and $dim\widetilde K_m =n$. We take some mapping
$\pi^m_{m-1}:\widetilde K_m \to \widetilde K_{m-1}$ generated by
$\omega_m$ inscribed into $\{ \varphi_{m-1}\wedge\eta\wedge
\gamma_m\}$ simplicial with respect to triangulation $K_m$ and
$K^*_{m-1}$.

Then, according to section b) of the remark 2, the mapping
$\pi^m_{m-1}f_m$ is descent of mapping $f_{m-1}$ with respect to
triangulation  $K^*_{m-1}$. Therefore
\begin{equation}
d\left( f_{m-1}, \pi^m_{m-1}f_m\right) < \frac1{2^m},\quad d\left( \pi_j^{m-1}f_{m-1}, \pi^m_{j}f_m\right) <
\frac1{2^m},\quad j<m-1.
\end{equation}

Continuing construction $n$-dimensional polyhedrons $\widetilde
K_i$ and mappings $\pi_i^{i+1},$ we obtain inverse sequence
 $\stackrel{.}{s} =\left\{
\widetilde K_i, \pi_i^{i+1}\right\},$ $i=1,2,\dots,$ satisfying
all conditions of theorem. By  $\widetilde s$ we denote limit of
inverse sequence  $s.$ Consider for each $i=1,2,3,\dots$ the
sequence of mappings
\begin{equation}
f_i,\, \pi_i^{i+1}f_{i+1},\, \pi_i^{i+2}f_{i+2},\, \dots
\end{equation}
of the space $X$ into polyhedron $\widetilde K_i.$ The proof of
that fact, which all later mappings of the sequence (2) are
obtained from  $f_i$ by descent with respect to triangulation
 $K_i$, similarly compact case of the space   $X$ (see [1, Chapter 5, \S 5,
 Freudenthal's Theorem]). According to second inequality of (1) we
  have
$$
d\left( \pi_i^{m-1} f_{m-1}, \pi_i^m f_m\right) < \frac1{2^m}
$$
Therefore for any point  $x\in X$ the sequence $\left\{ \pi_i^m
f_m(x)\right\},$ $m=i+1, i+2, \dots,$ is fundamental sequence. Since
the polyhedron $\widetilde K_i$ is complete metrizable, then the
sequence $\left\{ \pi_i^m f_m(x)\right\},$ $m=i+1,\dots,$ is
convergent at some point $g_i(x)\in \widetilde K_i.$ Sequence of
mappings $\left\{ \pi_i^m f_m\right\},$ $m=i+1, i+2,\dots,$ is
convergent to $g_i$ uniformly, therefore mapping $g_i: X\to
\widetilde K_i$ is continuous. Since all mappings  $ \pi_i^m f_m$
are obtained $f_i$ by descent with respect to triangulation $K_i,$
then the mapping $g_i$ also has this property (see [1, Chapter 4, \S
1, lemma 2]). Therefore mappings $g_i: X\to \widetilde K_i,$
$i=1,2,\dots,$ are canonical mappings with respect to covering
$\omega_i.$

Furthermore, according to preposition 3, mappings $g_i: X\to
\widetilde K_i,$ $i=1,2,\dots,$  are irreducible mappings with
respect to triangulation $K_i$(and, so, will be mappings on
$\widetilde K_i$). The relation $g_i=\pi_i^jg_j$ when $i<j$ is
checked by standard way (see [1, Chapter 5, \S 5]).

Since each mapping $g_i$ is  $\omega_i$-mapping of the space  $X$
into polyhedron  $\widetilde K_i,$ and system of open coverings
$\omega_i,$ $i=1, 2, \dots,$ of the space $X$ is refinement (see
[1, Chapter 1, \S 7, definition 10]) (since the covering
  $\omega_i$ is inscribed  in $\gamma_i$ ), then limit $g:X\to
\widetilde S \subseteq \prod\limits_{i=1}^\infty \widetilde K_i$
of mappings $g_i$ is (see [1, Chapter 6, \S 4, lemma 2]) imbedding
of the space $X$ into limit $\widetilde S$ of inverse sequence
$S.$ We prove, that $g$ there is mapping of the space $X$ on limit
$\widetilde S$ of the inverse sequence $S.$

We take some point  $y^0\in \widetilde S$ and assume $y^0=\left\{
y_i^0, i=1,2,\dots\right\}.$ Consider closed sets
$\Phi_i=g_i^{-1}y_i^0,$ $i=1,2,\dots $ in $X$. Since $g_i$ is
$\omega_i$-mapping, then $\Phi_i\subseteq O_{\alpha(i)}\in
\omega_i,$ $i=1,2,\dots$ . We prove, that $\Phi_{i+1}\subseteq
\Phi_i,$ $i=1, 2, \dots$ .

Since  $y_i^0=\pi_i^{i+1}y^0_{i+1},$ then $$ y^0_{i+1}\subseteq
\left( \pi_i^{i+1}\right)^{-1} y_i^0, \quad i=1,2,\dots\ .
\eqno(*)$$ Then from  inclusion  $(*)$  and the equality
$g_i=\pi_i^{i+1}g_{i+1}$ follows that
$$
g^{-1}_{i+1}y^0_{i+1}=\Phi_{i+1}\subseteq g^{-1}_{i+1}\left( \pi_i^{i+1}\right)^{-1} y_i^0=g_i^{-1}y_i^0 = \Phi_i,
\quad i=1,2,\dots \ .
$$
So, the system $\{ \Phi_i, i=1,2,\dots\}$ closed in $X$ sets
$\Phi_i,$ the sets which diameters tends to zero, is embedded. Then
from completeness of the space $X$ follows that intersection of the
sets  $\Phi_i$ nonempty and consists on one point. Suppose
$\bigcap\limits_{i=1}^\infty \Phi_i=\{ x^0\}.$ Since
$g_i\Phi_i=y_i^0$ and  $x^0\in \Phi_i,$ then $g_i(x^0)=y_i^0,$
$i=1,2,\dots$ . Consequently, $gx^0=y^0$ and therefore $y^0\in gX.$
Since $y^0$ is any point of the space   $\widetilde S,$ then from
here follows, that $g$ is (topological) mapping of the space $X$
onto limit $\widetilde S$ of the inverse  sequence  $S.$

Note that in identification of points $x\in X$ and $gx\in \widetilde
S$ projections $\pi_i:\widetilde S\to \widetilde K_i$ are identified
with irreducible with respect to triangulation $K_i$ mappings $g_i.$
\hfill $\Box$

This theorem is generalization of the Freudenthal's theorem [5].

{\bf Corollary 1.} Any $n$-dimensional metric  superparacompact
space $X$ is homeomorphic  to the everywhere dense subset of the
limit $\widetilde S$ of the inverse sequence $S=\left\{ \widetilde
K_i, \pi_i^{i+1}\right\},$ $i=1,2,\dots,$ from $n$-dimensional
polyhedron $\widetilde K_i,$ being bodies of standard triangulation
$K_i$ and decomposing into discrete sum of compact polyhedrons; in
addition projections $\pi_i^{i+1}$ are simplicial with respect to
$K_{i+1}$ and some triangulation $K_i^*$ of the polyhedron
$\widetilde K_i,$ being subdivision of the triangulation  $K_i.$
Each projection $\pi_i:X\to \widetilde K_i$ is irreducible with
respect to the triangulation  $K_i,$ $i=1,2,\dots$ .

\textbf{Proposition 4.} \emph{Any superparacompact complete with
respect to Cech   $(p-)$ space $X$ [6] is perfectly mapped into
Baire space  $B(\tau)$ of the weight $\tau$ (onto 0-dimensional in
the sense $dim$ metrizable space of the weight  $\leq\tau$).}

\textbf{Proof.} The space $X$ is perfectly mapped   (see [7,
theorem 2]) onto 0-dimensional in the sense $dim$ complete
metrizable (metrizable) space $X_0.$ Therefore $\omega
X_0\le\tau.$ Since any 0-dimensional in the sense  $dim$ complete
metrizable space of the weight $\le\tau$ is homeomorphic (see [8,
preposition  5.1]) closed subspace of generalized  Baire space
$B(\tau)$ of the weight $\tau$ and composition perfect mappings
are perfect, then hence follows, that the space $X$ is perfectly
mapped into Baire space $B(\tau)$ of the weight   $\tau$ (onto
0-dimensional in the sense $dim$ metrizable space of the weight
$\le \tau$).\hfill $\Box$

{\bf Corollary 2.} Any superparacompact complete metrizable space
$X$ of the weight  $\le\tau$ is perfectly mapped into Baire space
$B(\tau)$ of the weight  $\tau$.

\textbf{Theorem 2.} \emph{For metrizable space $X$ following
statements are equivalent: a) $X$  is superparacompact  complete
metrizable space of weight $\le\tau$; b) $X$ is perfectly mapping
into Baire space   $B(\tau)$ of the weight $\tau$; c) $X$ is closed
included into product $B(\tau)\times Q^\infty$ of Baire space
$B(\tau)$ of the weight $\tau$ on Hilbert cub  $Q^\infty$.}

\textbf{Proof.} If in the condition of the theorem
$\tau<\aleph_0,$ then all statements of the theorem are evidently.
Therefore we consider the case, when $\tau\ge\aleph_0.$

The statement b) implies from statement a) because of preposition
4.

The case b) $\Rightarrow$ c). Let $f$ be perfect mapping of the
space  $X$ into Baire space  $B(\tau)$ of the weight  $\tau$.
There exists (see [9, theorem 3]) such imbedding $g:X\to
B(\tau)\times Q^\infty,$ that $f=\pi \circ g,$ where $\pi$ is the
projection  $B(\tau)\times Q^\infty$ onto $B(\tau).$ Since the
mapping $f$ is perfect, and the space $B(\tau)\times Q^\infty$ is
Hausdorff space, then the mapping $g$ is perfect [4]. Thus, $g$ is
closed imbedding of the space $X$  into product $B(\tau)\times
Q^\infty$ of Baire space  $B(\tau)$ of the weight $\tau$ to
Hilbert cub $Q^\infty.$

Now we derive from statement c) the statement a). The product
$B(\tau)\times Q^\infty$ of Baire space $B(\tau)$ of the weight
 $\tau$ to  Hilbert cub $Q^\infty$ is superparacompact  (see
[3, corollary 1]). It is known, that the product $B(\tau)\times
Q^\infty$ is complete metrizable and $w\left(B(\tau)\times
Q^\infty\right)=\tau$. Then from monotonicity of complete metrizable
and superparacompact (see [3]) by closed subspaces follows, that the
space  $X$ is superparacompact and complete metrizable. Since
$w\left(B(\tau)\times Q^\infty\right)=\tau$, then $wX\le\tau.$\hfill
$\Box$

We note, that theorem 2 is extension of the theorem Morita  [10]
about universality of the product $B(\tau)\times Q^\infty$ in the
class of all strongly metrizable space of the weight $\le\tau.$

\textbf{Theorem 3.} \emph{For Hausdorff space  $X$ following
statement are equivalent: a) $X$ is superparacompact (complete)
metrizable space of the weight $\le\tau$ and $dimX\le n;$ b) $X$ is
closed imbedded into product  (Baire space  $B(\tau)$ of the weight
$\tau$) of 0-dimensional in the sense $dim$ of metrizable space of
the weight $\tau$ onto universal $n$-dimensional compact $\Phi^n.$}

\textbf{Proof.} By virtue of preposition 4, the space $X$ is
perfectly mapped (into Baire space $B(\tau)$ of the weight $\tau$)
onto 0-dimensional in the sense $dim$ metrizable space    $X_0$ of
the weight $\le\tau.$

Since the space $X$ is strongly metrizable, $dimX\le n$ and
$wX\le\tau,$ then by virtue of Nagata's theorem  (see [11]), the
space $X$ is topological mapped into product $B(\tau)\times
\Phi^n$ of generalized Baire space $B(\tau)$ of the weight  $\tau$
to universal   $n$-dimensional compact $\Phi^n.$ Then the space
 $X$ is homomorphic  (see [12, preposition 59, Chapter VI, \S
2]) to closed subspace of the product $\left(B(\tau)\times
B(\tau)\times \Phi^n\right) X_0\times B(\tau)\times \Phi^n$. Suppose
($B(\tau)=B(\tau)\times B(\tau)$) $R_\tau^0= X_0\times B(\tau).$ The
space $(B(\tau))R_\tau^0$ (complete) metrizable, ($wB(\tau)=\tau$)
$wR_\tau^0=\tau$ and 0-dimensional in the sense $dim$ [4].

We deduce from  statement b) the statement  a). The product
$\left( B(\tau)\times \Phi^n\right)$ $R_\tau^0\times \Phi^n$ is
superparacompact (see [3, corollary 1]). It is known [1], that the
product $\left( B(\tau)\times \Phi^n\right) R_\tau^0 \times
\Phi^n$ (complete) metrizable, $\left( dim\left( B(\tau)\times
\Phi^n\right)=n\right)$ $dim\left( R_\tau^0\times \Phi^n\right)=n$
and $\left( w\left( B(\tau)\times \Phi^n\right)=\tau\right)$
$w\left( R_\tau^0 \times \Phi^n\right)=\tau.$ then from
monotonicity of the superparacompact property   (see [3]),
complete metrizability of dimensionality $dim$ by closed subspaces
follows, that the space $X$ is superparacompact, (complete)
metrizable and $dimX\le n.$ Clearly, that and $wX\le\tau.$ \hfill
$\Box$

Theorem 3 is expansion of the Nagata's theorem  [11] about
embedding $n$-dimensional strongly metrizable space in
$B(\tau)\times \Phi^n$ to the case superparacomact.

\end{document}